\documentclass[a4paper,12pt,leqno]{article}
\usepackage{amsmath,amsthm,amssymb}
\usepackage[dvips]{graphicx}
\usepackage{indentfirst}
\numberwithin{equation}{section}
\allowdisplaybreaks

\newtheorem{thm}{Theorem}[section]
\newtheorem{lem}{Lemma}[section]
\newtheorem{prop}{Proposition}[section]

\newtheorem{definition}{Definition}[section]
\newtheorem{remark}{Remark}[section]
\newtheorem{example}{Example}[section]

\newtheorem{exercise}{Exercise}[section]
\newtheorem{report}{Report}[section]

\newcommand{\Proof}{\noindent {\bf Proof: }\ }

\newcommand{\bin}{\mathrm{Bin}}
\newcommand{\prob}{\mathrm{Prob}}
\newcommand{\kyori}{\mathrm{dist}}

\newcommand{\abslr}[1]{\lvert #1 \rvert}

\newcommand{\real}{\mathbb{R}}
\newcommand{\eqv}{\equiv}
\newcommand{\indicator}{\mathrm{1}}
\newcommand{\jouken}[1]{\mathrm{Condition \; {#1}}}
\newcommand{\btheta}{\bar{\theta}}

\newcommand{\barf}{\bar{f}}
\newcommand{\TrueParam}{T_{0}}
\newcommand{\Sni}{S'}

\title{Strong consistency of MLE for
finite uniform mixtures when the scale 
parameters are exponentially small}
\author{Kentaro Tanaka$^{1}$ and Akimichi Takemura$^{2}$}
\date{}

\begin{document}
\maketitle

\begin{center}
 {\scriptsize \textit{$^{1}$The Department of Industrial Engineering 
 and Management,  
 Tokyo Institute of Technology,  
 2-12-1 Ookayama, Meguro-ku,
 Tokyo 152-8552, JAPAN
 \\ 
 $^{2}$Department of Mathematical Informatics,  
 Graduate School of Information Science and Technology,  
 University of Tokyo, Bunkyo-ku, Tokyo 113-0033, Japan 
 }}
\end{center}

\begin{abstract}
We consider maximum likelihood estimation 
of finite mixture of uniform distributions.  We prove that maximum
likelihood estimator is strongly consistent, if the scale parameters
of the component uniform distributions are restricted from below by 
$\exp(-n^d)$, $0 < d < 1$, where $n$ is the sample size.
\end{abstract}


\section{Introduction}
\label{sec:intro}

\begin{figure}[ht]
 \begin{center}
  \input{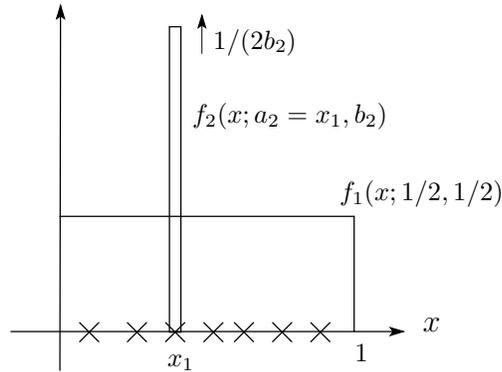} 
 \end{center}
 \caption{The likelihood tends to infinity as 
$b_{2} \rightarrow 0$ at $a_{2} = x_{1}$. }
 \label{unknown-d.tpic:fig}
\end{figure}

Consider a mixture of two uniform distributions
$$
 (1-\alpha) f_1(x;a_{1},b_{1}) + 
  \alpha f_2(x;a_{2},b_{2}) , 
$$
where 
$f_{m}(x;a_{m},b_{m})$, $m=1,2$,  are uniform densities 
with parameter $(a_m,b_m)$
on the half-open intervals $[a_{m} - b_{m}, a_{m} + b_{m})$ 
and $0 \leq \alpha \leq 1$.
For definiteness and convenience we use the half-open intervals in
this paper, although obviously the intervals can be open or closed.
By using half-open intervals, our densities are right continuous and
the version of the density is uniquely determined.
For simplicity suppose that $a_{1} = 1/2, b_{1} =1/2, \alpha =
\alpha_{0}$ are known and
the parameter space is
$$
  \{
  (a_{2}, b_{2})
  \mid
  0 \leq a_{2} -  b_{2} \; , \;
  a_{2} +  b_{2} \leq 1
  \}
$$
so that the support of the density is $[0,1)$.
Let $x_{1},\ldots,x_{n}$ denote 
a random sample of size $n \geq 2$ 
from the true density 
$(1-\alpha_{0}) f_1(x;1/2,1/2) + \alpha_{0} f_2(x;a_{2,0},b_{2,0})$. 
 If we set $a_{2} = x_{1}$, 
then likelihood tends to infinity as 
$b_{2} \rightarrow 0$ 
 (Figure~\ref{unknown-d.tpic:fig}). 
Hence the maximum likelihood estimator 
is not consistent. Actually it does not even
exist for each finite $n$.

When we restrict that $b_{2} \geq c$, where $c$ is a positive real
constant, then we can avoid the divergence of the likelihood and the
maximum likelihood estimator is strongly consistent provided that $b_{2,0} \geq
c$.  But there is a problem of how small we have to choose $c$ to
ensure $b_{2,0} \geq c$ 
since we do not know $b_{2,0}$.
An interesting question here is whether we
can decrease the bound $c=c_n$ to zero with the sample size $n$ and
yet guarantee the strong consistency of maximum likelihood
estimator.  If this is possible, the further question is how fast
$c_n$ can decrease to zero.  This question is similar to the (so far
open) problem stated in Hathaway(1985), which treats mixtures of
normal distributions with constraints imposed on the ratios of
variances.  See also a discussion in Section 3.8 of McLachlan and
Peel(2000).

\begin{figure}[ht]
 \begin{center}
  \includegraphics[width=13cm,clip]{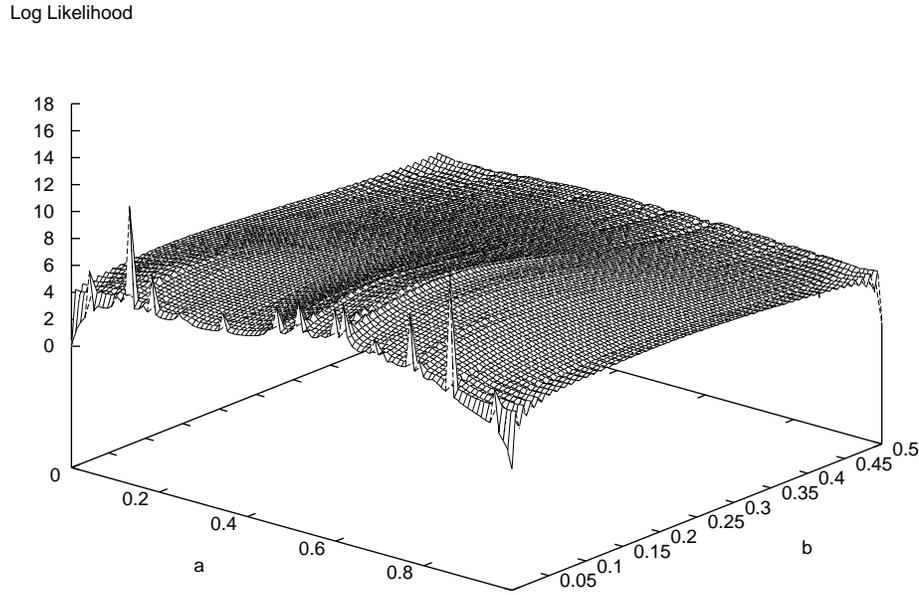}
 \end{center}
 \caption{An example of log likelihood function for $n=40$}
 \label{simlation:fig}
\end{figure}
Figure \ref{simlation:fig} depicts an example of
likelihood function. 
Random sample of size $n=40$ is generated from  
$0.6 \cdot f(x;0.5,0.5) + 0.4 \cdot f(x;0.6,0.2)$ and 
the model is $0.6 \cdot f(x;0.5,0.5) + 0.4 \cdot f(x;a,b)$.
Despite the limited resolution in Figure \ref{simlation:fig}
, there are actually $n=40$ peaks of the
likelihood function as $b \downarrow 0$.
We see that although the
likelihood function diverges to infinity at these peaks, 
the divergence takes place only for very small $b$ and the likelihood
function is well-behaved for most of the ranges of $b$.
This suggests that
the bound $c_n$ can decrease to zero fairly quickly 
while maintaining the consistency of maximum likelihood estimator.
In fact we
prove that $c_n$ can decrease exponentially fast to zero for the
mixture of $M$ uniform distributions.  More precisely we prove that 
maximum likelihood estimator is strongly consistent if
$c_n=\exp(-n^d)$, $0 < d < 1$.

The organization of the paper is as follows. 
In Section \ref{sec:preliminaries} we summarize 
some preliminary results. 
In Section \ref{sec:main} we state our main result 
in Theorem \ref{thm:main}. 
Proof of Theorem \ref{thm:main} is given in Appendix \ref{sec:proof}. 
In Section \ref{sec:discussion} we
give a simulation result and some discussions. 

\section{Preliminaries on identifiability of 
mixture distributions and strong consistency}
\label{sec:preliminaries}

In this section, we consider the identifiability and strong
consistency of finite mixtures.  The properties of finite mixtures
treated in this section concerns general finite mixture distributions.

A mixture of $M$ densities with parameter 
$\theta=(\alpha_1, \eta_1,\ldots,\alpha_M, \eta_M)$ is defined by 
$$
  f(x;\theta) \eqv \sum_{m=1}^{M} \alpha_{m} f_{m}(x;\eta_{m}) , 
$$
where $\alpha_{m}$, $m=1,\ldots,M$, called the mixing weights, 
are nonnegative real numbers that sum to one and 
$f_{m}(x;\eta_{m})$ are densities  
with parameter $\eta_{m}$.
$f_{m}(x;\eta_{m})$ are called the components 
of the mixture.
Let $\Theta$ denote the parameter space. 

In general, identifiability of a parametric family of densities  
is defined as follows. Note that in this paper a version of the
density is uniquely determined by the right continuity.
\begin{definition} 
 $($identifiability of a parametric family of densities$)$ \\
A parametric family of densities  
$\{f(x;\theta) \mid \theta \in \Theta\}$ 
is identifiable if 
different values of parameter designate different densities; 
that is 
\begin{equation} 
  f(x;\theta) = f(x; \theta')
  \quad \forall x, 
  \nonumber 
\end{equation}
implies $\theta = \theta'$. 
\label{identifiability:def}
\end{definition}
If a parametric family of densities is not identifiable, 
then it is said to be unidentifiable. 

In mixture case, 
when all components $f_{m}(x;\eta_{m})\; ,\; m=1,\ldots ,M$ 
belong to the same parametric family, 
then $f(x;\theta)$ is invariant 
under the permutations of the component labels. 
Because of this trivial unidentifiability, 
the definition of identifiability for the mixture densities 
can be weakened as described in 
Teicher(1960), Yakowitz and Spragins(1968), McLachlan and Peel(2000) 
and so on, 
so that 
$
\sum_{m=1}^{M}\alpha_{m}f_{m}(x;\eta_{m}) 
 = \sum_{m'=1}^{M'}\alpha_{m'}'f_{m'}(x;\eta_{m'}')
$
 implies $M = M'$ and for each $m$ there exists some $m'$ 
such that $\alpha_{m} = \alpha_{m'}$ 
and $\eta_{m} = \eta_{m'}'$. 
But, even under such a weakened definition, 
mixtures of density functions still have unidentifiability. 
For example, if $\alpha_{1} = 0$, then for all
parameters which differ only in $\eta_{1}$, we have the same
density. 
We also discuss examples of non-trivial unidentifiability of mixtures 
after theorem~\ref{thm:main} below. 
In any way, mixture model is unidentifiable. 

In unidentifiable case, true model may consist of two
or more points in the parameter space.  Therefore we have to carefully
define strong consistency of estimator $\hat\theta_n$, 
because we should define $\hat\theta_n$ to be consistent if $\hat\theta_n$
falls in arbitrary small neighborhood of the set of points
designating the true model as $n \rightarrow\infty$.

The following definition is
essentially the same as Redner's(1981).  We suppose that the parameter
space $\Theta$ is a subset of Euclidean space and
$\kyori(\theta,\theta')$ denotes the Euclidean distance between
$\theta,\theta' \in \Theta$.
\begin{definition} 
 $($strongly consistent estimator$)$ \\
 Let $\TrueParam$ denote the set of true parameters
 \begin{eqnarray}
  \TrueParam \eqv \{\theta \in \Theta \mid f(x;\theta) = f(x;\theta_0)
   \quad \forall x 
   \} , 
   \nonumber
 \end{eqnarray}
 where $\theta_{0}$ is one of parameters 
 designating the true distribution.
 An estimator $\hat{\theta}_n$ is strongly consistent if
  \begin{eqnarray}
  \prob
   \left(
    \lim_{n \rightarrow \infty}
    \inf_{\theta \in \TrueParam}\kyori(\hat{\theta}_n, \theta) = 0
   \right) = 1.
   \nonumber
 \end{eqnarray}
\end{definition}

In this paper two notations $\prob(A)=1$ and $A, a.e.$
($A$ holds almost everywhere), will be used interchangeably.
The index $_{0}$ to the parameter 
always denotes the true parameter.

In finite mixture case,
regularity conditions for 
strong consistency of maximum likelihood estimator
are given in Redner(1981).
When the components of the mixture are the densities 
of continuous distributions
and the parameter space is Euclidean, 
the conditions become as follows. 
Let $\Gamma$ denote a subset of the parameter space.

\begin{description}
 \item[$\jouken{1}$.] $\Gamma$ is a compact subset of Euclidean space.
\end{description}

For $\theta \in \Gamma$ and any positive real number $r$, let 
\begin{eqnarray}
f(x;\theta,r) & = & \sup_{\kyori(\theta',\theta) \leq r}f(x;\theta'), 
\nonumber \\
f^{\ast}(x;\theta,r) & = & \max(1,f(x;\theta,r))
\; . 
\nonumber
\end{eqnarray}
\begin{description}

 \item[$\jouken{2}$.] 
            For each $\theta \in \Gamma$ and sufficiently small $r$, 
            $f(x;\theta,r)$ is measurable and 
            \begin{equation}
             \int 
                   \log(f^{\ast}(x;\theta,r))
             f(x;\theta_0)dx < \infty
             \; .
            \end{equation}
\end{description}
\begin{description}
 \item[$\jouken{3}$.] 
            If
            $\lim_{n \rightarrow \infty}\theta_{n} = \theta$, 
            then 
            $\lim_{n \rightarrow \infty}
            f(x;\theta_{n}) = f(x;\theta)$
            except on a set which is a null set and 
            does not depend on the sequence 
            $\{\theta_{n}\}_{n=1}^{\infty}$.
\end{description}
\begin{description}
 \item[$\jouken{4}$.]
            \begin{equation}
             \int \abslr{
              \log{f(x;\theta_0)}
              }
              f(x;\theta_0)dx < \infty .
            \end{equation}
\end{description}
The following two theorems have been proved by 
Wald(1949), Redner(1981).
\begin{thm}
 $\mathrm{(Wald(1949),Redner(1981))}$ 
 Suppose that {\rm Conditions  1, 2, 3 and 4} are satisfied. 
 Let $S$ be any closed subset of $\Gamma$ 
 not intersecting $\TrueParam$. Then
\begin{equation}
 \prob
  \left(\lim_{n \rightarrow \infty}
   \frac{ \sup_{\theta \in S}
   f(x_1;\theta)\times \cdots \times f(x_n;\theta)} 
   {f(x_1;\theta_0)\times \cdots \times f(x_n;\theta_0)}
   = 0
  \right)
  = 1
  \; .
  \label{WaldAndRedner:eq}
\end{equation}
\label{WaldAndRedner:thm}
\end{thm}
\begin{thm}
 $\mathrm{(Wald(1949),Redner(1981))}$ 
 Let $\tilde{\theta}_{n}$ be any function of the 
 observations $x_{1},\ldots,x_{n}$ such that
 \begin{eqnarray}
  \forall n, \quad
  \prod_{i=1}^{n}
   \frac{f(x_{i};\tilde{\theta}_{n})}{f(x_{i};\theta_{0})}
   \geq \delta > 0 , 
   \nonumber 
 \end{eqnarray}
 then 
 $\prob( 
 \lim_{n \rightarrow \infty} 
 \inf_{\theta \in T_{0}}{\kyori(\tilde{\theta}_{n},\theta)}
 ) = 1 .$
\label{Wald:thm}
\end{thm}
If Conditions 1, 2, 3 and 4  are satisfied, 
then it is readily verified 
by theorems \ref{WaldAndRedner:thm} and \ref{Wald:thm} 
that maximum likelihood  estimator restricted to $\Gamma$ is 
strongly consistent. 

We also state Okamoto's inequality, which will be used in our proof in
Appendix \ref{sec:proof}.

\begin{thm}
\label{thm:Okamoto}
$\mathrm{(Okamoto(1958))}$
 Let $Z$ be a random variable following a binomial distribution $\bin(n, p)$.
 Then for $\delta > 0$  
 \begin{eqnarray}
  \prob\left(\frac{Z}{n} - p \geq \delta \right) <
  \exp{(-2n\delta^2)}.
  \label{okamoto1:ineq}
 \end{eqnarray}
\end{thm}

\section{Main result}
\label{sec:main}

Here, we generalize the problem stated in introduction 
to the problem of mixture of $M$ uniform distributions and 
then state our main theorem.

A mixture of $M$ uniform densities 
with parameter 
$\theta$
is defined by 
\begin{eqnarray}
f(x;\theta) \eqv \sum_{m=1}^{M} \alpha_{m} f_{m}(x;\eta_{m}), 
\nonumber
\end{eqnarray}
where $f_{m}(x;\eta_{m}) \eqv f_{m}(x;a_{m},b_{m})$, $m=1,\ldots,M$,
are uniform densities 
with parameter $\eta_{m} = (a_m,b_m)$
on half-open intervals $[a_{m} - b_{m}, a_{m} + b_{m})$ 
and $\alpha_{m}$ are mixing weights.
The parameter space $\Theta \subset \real^{3M}$ is defined by
$$
 \Theta
  \eqv
  \{
  (\alpha_1,a_1,b_1,
  \ldots,\alpha_M,a_{M}, b_{M})
  \mid
  0 \leq \alpha_{1},\ldots,\alpha_{M} \leq 1 \; , \;
  \sum_{m=1}^{M}\alpha_{m} = 1 \; , \; 
  b_{1},\ldots,b_{M} > 0
  \} \; .
$$

Let 
$\theta_{0} \eqv (\alpha_{0,1},a_{0,1},b_{0,1},
  \ldots,\alpha_{0,M},a_{0,M}, b_{0,M}) $ 
be the true parameter and let
$$
f(x ; \theta_0)= \sum_{m=1}^{M} \alpha_{0,m}\; f_{m}(x; a_{0,m}, b_{0,m})
$$
be the true density.  
Denote the minimum and the maximum of the support of $f(x ; \theta_0)$
by
\begin{eqnarray*}
  L_{\min} & = & \min(a_{0,1} - b_{0,1},\ldots,a_{0,M} - b_{0,M}) , \\
  L_{\max} & = & \max(a_{0,1} + b_{0,1},\ldots,a_{0,M} + b_{0,M}) ,
\end{eqnarray*}
and let
$$
  L  = L_{\max} - L_{\min}.
$$

Let $\Theta_{c}$ be a constrained parameter space
\begin{eqnarray}
 \Theta_{c}
  \eqv
  \{
  \theta \in \Theta
  \mid
  b_{m} \geq c > 0
  \; , \;
  m=1,\ldots,M
  \} , 
  \nonumber
\end{eqnarray}
where $c$ is a positive real constant.
We can easily see that 
Conditions 1, 2, 3 and 4 are satisfied with $\Theta_c$.  Therefore
if $\theta_{0} \in \Theta_{c}$, 
then maximum likelihood  estimator restricted to $\Theta_{c}$
is strongly consistent (Redner(1981)).
But there is a problem of 
how small $c$ must be to ensure 
$\theta_{0} \in \Theta_{c}$ as discussed in section
\ref{sec:intro}.

Since the support of uniform density is compact,
the following lemma holds.

\begin{lem}
\label{boundedParameterSpace:lem}
 For any parameter 
 $\theta = (\alpha_{1},a_{1},b_{1},\ldots,
 \alpha_{M},a_{M},b_{M}) \in \Theta$,
 there exists a parameter 
 ${\theta}' = (\alpha_{1},a_{1}',b_{1}',\ldots,
 \alpha_{M},a_{M}',b_{M}') \in \Theta$
 satisfying  
$$
L_{\min} \leq 
   a_{1}',\ldots,a_{M}'
   \leq L_{\max}, \quad
  0 < b_{1}',\ldots,b_{M}'
   \leq L
$$
such that
$$
   \sum_{m=1}^{M} \alpha_m f_{m}(x; a_{m}', b_{m}')
   \geq
   \sum_{m=1}^{M} \alpha_m f_{m}(x; a_{m}, b_{m}), 
\quad \forall x\in [L_{\min}, L_{\max}), 
$$
where equality does not hold 
if there exists $\alpha_m > 0$ such that
$a_{m} \not\in 
[L_{\min}, L_{\max}]$ 
or $b_{m} > L$.
\end{lem}

By $\mathrm{lemma~\ref{boundedParameterSpace:lem}}$,
maximum likelihood  estimator is 
restricted to a bounded set in $\Theta \subset\real^{3M}$.

Let $\{c_n\}_{n=0}^{\infty}$ be a
monotone decreasing sequence of positive real numbers 
converging to zero 
and define $\Theta_{n}$ by 
$$
 \Theta_{n}
   \eqv 
  \{
  \theta \in \Theta
  \mid
  0 < c_{n} \leq b_{m} , \
  m=1,\ldots,M
  \}
  \; .
$$

We are now ready to state our main theorem.

\begin{thm} 
 \label{thm:main}
 Suppose that the true model $f(x;\theta_0)$ can not be represented  by 
 any model consisting of less than $M$ components.
 Let $c_{0}> 0 $ and $0 < d < 1$.
 If $c_{n} = c_{0} \exp{(-n^d)} \leq b_{m}$ for all $b_{m}$,   
 then maximum likelihood estimator (which is restricted to $\Theta_{n}$)
 is strongly consistent. 
\end{thm}
Proof of this theorem is given in Appendix \ref{sec:proof}.

Note that under the assumption of theorem~\ref{thm:main} the strong
consistency holds even if the true model is unidentifiable in a
non-trivial way.  We illustrate the assumption of
theorem~\ref{thm:main} by examples of two-component models.
If the true model is 
$\alpha U(x;0,\alpha) + (1-\alpha) U(x;\alpha,1)$ 
(see Titterington et.\ al.\ (1985) pp.\ 36)
which is unidentifiable 
and can be represented by one component model, 
then the assumption of theorem Theorem \ref{thm:main} is not satisfied.
But if the true model is represented by 
$\frac{1}{3}U(x;-1,1) + \frac{2}{3}U(x;-2, 2)$
(see Everitt and Hand(1981) pp.\ 5), 
which is unidentifiable because 
$\frac{1}{2}U(x;-2,1) + \frac{1}{2}U(x;-1,2)$ 
represents the same distribution,  
then the assumption of theorem Theorem \ref{thm:main}
 is satisfied, because 
it  can not be represented by one component model. 

Next proposition states that 
the rate of $c_{n} = \exp(-n^d)$, $d < 1$, 
obtained in theorem \ref{thm:main}
is almost the lower bound of the order of $c_n$ which maintains
the consistency.

\begin{prop}
If $c_{n}$ decreases faster than $\exp(-n)$, i.e., $e^n c_n
\rightarrow 0$, 
then the consistency of maximum likelihood estimator 
restricted to 
$\Theta_{n}$ 
fails.
\end{prop}
\Proof
By the strong law of large numbers, 
mean log likelihood of true model 
\\
$\frac{1}{n}\log{\sum_{i=1}^{n}f(x_{i};\theta_{0})}$
converges to 
$ E[\log{f(x;\theta_{0})}] < \infty $
almost everywhere.
Assume that  $c_{n}$ decrease faster than $\exp(-n)$.
Take  $a_{1} = x_{1}, b_{1} = c_{n}$.  Fix $\alpha_{1} > 0$ and fix other
parameters
$(\alpha_2, \eta_2,\ldots,\alpha_M, \eta_M)$ such that
$\frac{1}{n}  \sum_{i=2}^{n}
    \log{
    \{
     \sum_{m=2}^{M} \alpha_{m}f_{m}(x_{i};\eta_{m})
    \}
    }
$ converges to a finite limit almost everywhere.
Then 
\begin{eqnarray}
 \lefteqn{
 \frac{1}{n}
  \sum_{i=1}^{n}
  \log{\sum_{m=1}^{M} f_{m}(x_{i};\eta_{m})}
  }& & 
  \nonumber \\ 
 & \geq & 
  \frac{1}{n}
   \log{
   \left\{\alpha_{1}f_{1}(x_{1};a_{1}=x_{1}, b_{1}=c_{n})\right\}
   }
   + 
   \frac{1}{n}\sum_{i=2}^{n}
    \log{
    \left\{
     \sum_{m=1}^{M} \alpha_{m}f_{m}(x_{i};\eta_{m})
    \right\}
    }
  \nonumber \\
&  \ge &
  \frac{1}{n}
  \log{\left\{
  \frac{\alpha_{1}}{2c_{n}}
  \right\}
  }
   + \frac{1}{n}
   \sum_{i=2}^{n}
   \log{
   \left\{
    \sum_{m=2}^{M} \alpha_{m}f_{m}(x_{i};\eta_{m})
   \right\}
   }
   \rightarrow \infty.
   \nonumber 
\end{eqnarray}
Therefore mean log likelihood of the true model is 
dominated by that of other models and  
consistency of maximum likelihood estimator fails.
\qed

\section{Some discussions}
\label{sec:discussion}

As stated above in Section \ref{sec:intro}, 
the failure of consistency 
of maximum likelihood estimator 
is caused by the divergence of the likelihood 
of the model, where  some scale parameters go to zero. 
Therefore in our setting it is of interest to investigate
the behavior of the likelihood of the models 
on the boundary ($b_{m} = c_{n}$) 
of the restricted parameter space $\Theta_{n}$. 
We report a simulation result for the case that 
the true model is 
$0.6 \cdot f(x;0.5,0.5) + 0.4 \cdot f(x;0.6,0.2)$ and 
a competing model is $0.6 \cdot f(x;0.5,0.5) + 0.4 \cdot f(x;a,b=c_{n})$ 
which is on the boundary ($b = c_{n}$) of the restricted parameter space, 
where $c_{n} = \exp(n^{-0.93})$.
\begin{table}[htbp]
\caption{log likelihood of the true model and 
that of a competing model}
\label{sim:tab}
\begin{center}
\begin{tabular}{|c|c|c|} \hline
 sample size $n$ & log likelihood (true)
 & log likelihood ($b= c_{n}$)\\ \hline
 10 & 0.7767 & 2.305 \\ \hline
 50 & 9.769 & 11.38 \\ \hline
 100 & 15.61 & 20.26 \\ \hline
 500 & 56.49 & 67.11 \\ \hline
 1000 & 117.9 & 104.7 \\ \hline
 5000 & 582.6 & 199.3 \\ \hline
\end{tabular}
\end{center}
\end{table}
The second column of Table \ref{sim:tab} shows the log likelihood at
$\hat\theta_n =\theta_0$.  The third column shows 
the log likelihood maximized with respect to $a\in [0,1]$ (but $b$ is taken to
be $c_n$).
In the competing model,  with probability tending to 1,
the length of the interval $2c_{n}$ 
is shorter than 
the minimum of the distance between realized values.
Therefore with probability tending to 1 the support of $f(x;a,b=c_{n})$ 
does not contain two or more realized values 
for all $a \in [0,1]$.
Therefore the maximum of the likelihood is usually achieved 
when the support of $f(x;a,b=c_{n})$ contains just one realized value.
Then $f(x;a,b=c_{n}) = 0.6 + 0.4/(2c_{n})$ on one particular realization
and $f(x;a,b=c_{n}) = 0.6$ on the other $n-1$ realized values. 
In this case 
the maximum of the log likelihood in competing model is  given by
$\log{\{0.6 + 0.4/(2c_{n})\}} + (n-1)\log\{0.6\}$.
The result in Table \ref{sim:tab} is based on 
one replication for each sample size. 
If we repeat the simulations, 
the results are similar.
Therefore 
the result in Table \ref{sim:tab} indicates that 
the log likelihood of the true model gets larger than 
that of the competing models with $b=c_n$  as the sample size $n$ increases. 
This simulation result is consistent with Theorem~\ref{thm:main}.

We expect that our result can be extended to other finite mixture
cases, especially for densities which are Lipschitz continuous when
the scale parameters are fixed. 
On the other hand, in Theorem~\ref{thm:main}, 
it might be difficult to weaken the assumption 
that there is no representation of the true model 
with less than $M$ components. 
The problem studied in this paper is
similar to the question stated in Hathaway(1985) which treats the
normal mixtures and the constraint is imposed on the ratios of
variances.  Methods used in this paper may be useful to solve the
question.

\appendix
\section{Appendix : Proof of the strong consistency}
\label{sec:proof}

Here we present a proof of Theorem \ref{thm:main}. 
Note that it is sufficient to prove Theorem \ref{thm:main} 
for $d$ arbitrarily close to 1. 
Therefore we assume $d > 1/4$ hereafter.

The whole proof is
long and we divide it into smaller steps. Intermediate results will be
given in a series of lemmas.

Define 
\begin{eqnarray}
  \Theta_{n}'
  & \eqv & 
  \{
  \theta \in \Theta_{n}
 \mid
  L_{\min} \leq \forall a_{m} \leq L_{\max}
  \; , \; 
  c_{n} \leq \forall b_{m} \leq L 
  \; , \; 
  c_{n} \leq \exists b_{m} \leq c_{0}
   \}, 
  \nonumber \\ 
 \Gamma_{0}
  & \eqv &  
  \{
  \theta \in \Theta
  \mid
  L_{\min} \leq a_{m} \leq L_{\max}
  \; , \; 
  c_{0} \leq b_{m} \leq L
  \; , \;
  m=1,\ldots,M
  \}
  \; .
  \nonumber
\end{eqnarray}
Because $\{c_n\}$ is decreasing
to zero, by replacing $c_0$ by some $c_n$ if necessary, 
we can assume without loss
of generality that $\TrueParam \subset \Gamma_0$.

In view of Theorems~\ref{WaldAndRedner:thm},~\ref{Wald:thm},
for the strong consistency of MLE on $\Theta_n$, 
by Lemma \ref{boundedParameterSpace:lem}, 
it suffices to prove that 
$$
 \lim_{n \rightarrow \infty}
  \frac{
   \sup_{\theta \in \Sni \cup \Theta_{n}'}
   \prod_{i=1}^{n} f(x_{i};\theta)
   }
   { \prod_{i=1}^{n} f(x_{i};\theta_{0}) } = 0, \quad a.e. 
$$
for all closed $\Sni \subset \Gamma_0$ not intersecting $\TrueParam$. 
Note that for all $\Sni$ and $\{x_{i}\}_{i=1}^{n}$, 
\begin{eqnarray}
 {
  \sup_{\theta \in \Sni \cup \Theta_{n}'}
  \prod_{i=1}^{n} f(x_{i};\theta)
  }
  = 
  \max
  \left\{
  {
  \sup_{\theta \in \Sni}
  \prod_{i=1}^{n} f(x_{i};\theta)
  }
  \; , \; 
  {
  \sup_{\theta \in \Theta_{n}'}
  \prod_{i=1}^{n} f(x_{i};\theta)
  }
  \right\}
  \; . \;
  \nonumber
\end{eqnarray}
Furthermore
equation (\ref{WaldAndRedner:eq}) 
with $S$ 
replaced by $\Sni$
holds by Theorem \ref{WaldAndRedner:thm}. 
This implies that it suffices to prove 
equation (\ref{WaldAndRedner:eq}) 
with $S$ replaced by $\Theta_{n}'$.

Note that in the argument above the supremum of the likelihood
function over $\Sni \cup \Theta_{n}'$ is considered separately for $\Sni$
and $\Theta_n'$.  $\Sni$ and $\Theta_n'$ form a covering of $\Sni\cup
\Theta_n'$.  In our proof, we consider finer and finer finite coverings of
$\Theta_n'$.  As above, it suffices to prove that the ratio of the
supremum of the likelihood over each member of the covering to the
likelihood at $\theta_0$ converges to zero almost everywhere.

Let $\theta \in \Theta_{n}'$.
Let $K \eqv K(\theta) \geq 1$ be the number of components
which satisfy ${b_{m}} \leq c_{0}$.  
Without loss of generality,
we can set 
$b_{1} \leq b_{2} \leq \cdots \leq b_{K}
\leq c_{0}  <  b_{K+1} \leq \cdots \leq b_{M}$.
Let $\Theta_{n,K}'$ be 
\begin{eqnarray}
 \Theta_{n,K}' 
  \eqv
  \{
  \theta \in \Theta_{n}'
  \mid
  b_{1} \leq b_{2} \leq \cdots \leq b_{K}
  \leq c_{0} < b_{K+1} \leq \cdots \leq b_{M}
  \}
  \; . 
  \nonumber 
\end{eqnarray}
Our first covering of $\Theta_n'$ is given by 
$$
\Theta_n' = \bigcup_{K=1}^M \Theta_{n,K}' .
$$
As above, it suffices to prove 
equation (\ref{WaldAndRedner:eq}) 
with $S$ replaced by $\Theta_{n,K}'$. 
We fix $K$ from now on.
Define $\bar{\Theta}_{K}$ by
\begin{eqnarray}
  \bar{\Theta}_{K} & \eqv &
  \{
  (\alpha_{K+1},a_{K+1},b_{K+1},\ldots,\alpha_{M},a_{M},b_{M})
  \in \real^{3(M-K)}
  \mid \sum_{m=K+1}^{M} \alpha_{m} \leq 1 \; , \;
  \alpha_{m} \geq 0 \; ,\; 
  \nonumber \\
 & & \hspace{3.8cm}
  L_{\min} \leq a_{m} \leq L_{\max} \; ,\;
  c_{0} \leq b_{m} \leq L \; ,\;
  m = K+1,\ldots,M
  \}
  \nonumber 
\end{eqnarray}
and for $\btheta \in \bar{\Theta}_{K}$, define 
\begin{eqnarray}
 \barf(x;\btheta) & \eqv & \sum_{m=K+1}^{M} \alpha_{m} f_{m}(x;\eta_{m})
  \; , \;
  \nonumber \\
 \barf(x;\btheta, \rho) & \eqv & 
  \sup_{\kyori(\btheta,\btheta') \leq \rho}
  \barf(x;\btheta')
  \; .
  \nonumber
\end{eqnarray}
Note that $\barf(x;\btheta)$ is a subprobability measure.
\begin{lem}
Let $B(\btheta,\rho(\btheta) )$ denote the open ball
with center $\btheta$ and radius $\rho(\btheta)$.
Then $\bar{\Theta}_{K}$ can be covered by 
a finite number of balls 
$B(\btheta^{(1)}, \rho(\btheta^{(1)})),\ldots,
B(\btheta^{(S)}, \rho(\btheta^{(S)}))$
such that 
\begin{equation}
 E_{0}[\log{\barf(x;\btheta^{(s)},\rho(\btheta^{(s)}))}]
  <
  E_{0}[\log{f(x;\theta_0)}]
  \; , \quad s=1,\ldots,S,
  \label{KLinequality:SUPbarf:wald:eq}
\end{equation}
where $E_{0}[\cdot]$ denotes the expectation under $\theta_{0}$.
\label{WaldThm1Type:lem}
\end{lem}

\Proof 
The proof is the same as in Wald (1949).
For all $\btheta \in \bar{\Theta}_{K}$, 
there exists a positive real number $\rho(\btheta)$ 
which satisfies 
$$
 E_{0}[\log{\barf(x;\btheta,\rho(\btheta))}] <
E_{0}[\log{f(x;\theta_0)]}. 
$$
Since 
$\bar{\Theta}_{K} \subset 
\bigcup_{\btheta}B(\btheta, \rho(\btheta))$
and $\bar{\Theta}_{K}$ is compact, 
there exists a finite number of balls 
$B(\btheta^{(1)}, \rho(\btheta^{(1)})),\ldots,
B(\btheta^{(S)}, \rho(\btheta^{(S)}))$
which cover $\bar{\Theta}_{K}$. 
\qed 

\bigskip
Define  
$$
 \Theta_{n,K,s}'
 \eqv 
  \{
  \theta \in \Theta_{n,K}'
  \mid
  (\alpha_{K+1},a_{K+1},b_{K+1},\ldots,\alpha_{M},a_{M},b_{M})
  \in B(\btheta^{(s)},\rho(\btheta^{(s)}))
  \}. 
$$
We now cover $\Theta_{n,K}'$ by 
$\Theta_{n,K,1}', \ldots, \Theta_{n,K,S}'$ :
$$
\Theta_{n,K}' = \bigcup_{s=1}^S \Theta_{n,K,s}' \   .
$$
Again it suffices to prove
that for each $s$, $s=1,\ldots,S$,
\begin{equation}
   \lim_{n \rightarrow \infty}
   \frac{
   \sup_{\theta \in \Theta_{n,K,s}'}
   \prod_{i=1}^{n} f(x_{i};\theta)
   }
   { \prod_{i=1}^{n} f(x_{i};\theta_{0}) }
   = 0, \quad a.e.
  \label{goal3:LikelihoodRatio:eq}
\end{equation}
We fix $s$ in addition to $K$ from now on.

Because 
$$
\lim_{n\rightarrow\infty} \frac{1}{n} \sum_{i=1}^{n} 
\log{f(x_{i};\theta_0)} =
E_0[\log f(x;\theta_0)], \quad a.e.
$$
(\ref{goal3:LikelihoodRatio:eq}) is implied by
\begin{eqnarray}
\limsup_{n \rightarrow \infty} \frac{1}{n} \sup_{\theta \in \Theta_{n,K,s}'}
   \sum_{i=1}^{n} \log{f(x_{i};\theta)}
  < 
   E_0[\log f(x;\theta_0)], 
  \quad a.e.  
  \label{goal:eq}
\end{eqnarray}
Therefore it suffices to prove (\ref{goal:eq}), which is a new
intermediate goal of our proof hereafter.

Choose $G$, $0 < G < 1$, such that
\begin{eqnarray}
 \lambda
  \eqv
  E_{0}[\log{f(x;\theta_0)}]
  - 
  E_{0}[\log{
  \{ \barf(x;\btheta^{(s)},\rho(\btheta^{(s)})) + G \}
  }] > 0
  \; . 
  \label{permissible_quantity:eq}
\end{eqnarray}
Let $u \eqv \max_{x} f(x;\theta_0)$.  Because $\{c_n\}$ is decreasing
to zero, by replacing $c_0$ by some $c_n$ if necessary, 
we can again assume without loss
of generality that $c_{0}$ is small enough to satisfy
\begin{eqnarray}
 2 c_0 &<& e^{-1} , \nonumber \\
 3M \cdot u \cdot 2c_{0} \cdot (- \log{G})
  & < &
  \frac{\lambda}{4} ,
  \label{1:c_0:condition:proof:eq}
  \\
 2{M} \cdot u \cdot 2c_{0} \cdot \log{\frac{1}{2c_{0}}}
  & < &
  \frac{\lambda}{12}
  \; .
  \label{2:c_0:condition:proof:eq}
\end{eqnarray}
Although $G$
depends on $c_0$, it can
be shown that $G$ and $c_0$ can be chosen small enough to satisfy 
these inequalities.
We now prove the following lemma.
\begin{lem}
Let $J(\theta)$ denote the support of 
$\sum_{m=1}^{K} \alpha_{m}f_{m}(x;\eta_{m})$
and let $R_n(V)$ denote the number of observations
which belong to a set $V \subset \real$.
Then for $\theta \in \Theta_{n,K,s}'$
\begin{eqnarray}
  \frac{1}{n}\sum_{i=1}^{n}\log{f(x_{i};\theta)}
  & \leq & 
  \frac{1}{n}\sum_{i=1}^{n}
  \log{
  \left\{
     \barf(x_{i};\btheta^{(s)},\rho(\btheta^{(s)})) + G
  \right\}
  }
  \\
  & & +
   \frac{1}{n}\sum_{x_{i} \in J(\theta)}
   \log{f(x_{i};\theta)}
   + 
   \frac{1}{n}R_n(J(\theta))\cdot (-\log{G})
   \; . \nonumber
   \label{loglikelihood:lem:eq}
\end{eqnarray}
\end{lem}

\Proof For $x\not\in J(\theta)$, 
$f(x;\theta)=\sum_{m=K+1}^M \alpha_m f_m(x;\eta_m)$.  Therefore
\begin{eqnarray}
\frac{1}{n}\sum_{i=1}^{n}\log{f(x_{i};\theta)}
& =& \frac{1}{n} \sum_{x_i \in J(\theta)} \log{f(x_{i};\theta)}
    + \frac{1}{n} \sum_{x_i \not \in J(\theta)} 
       \log\left\{\sum_{m=K+1}^M \alpha_m f_m(x_{i};\eta_m)\right\}
  \nonumber \\
&\leq &
  \frac{1}{n}\sum_{i=1}^{n}
  \log{
  \left\{
   \sum_{m=K+1}^{M} \alpha_{m} f_{m}(x_{i};\eta_{m}) + G
  \right\}
  }
  \nonumber \\
&& \qquad  +
   \frac{1}{n}\sum_{x_{i} \in J(\theta)}
   \left[
    \log
    f(x_i ; \theta)
    -
    \log
    \left\{
     \sum_{m=K+1}^{M} \alpha_{m} f_{m}(x_{i};\eta_{m}) + G 
    \right\}
   \right]
   \nonumber \\
&\leq &
  \frac{1}{n}\sum_{i=1}^{n}
  \log{
  \left\{
     \barf(x_{i};\btheta^{(s)},\rho(\btheta^{(s)})) + G
  \right\}
  }
  \nonumber \\
&& \qquad +
   \frac{1}{n}\sum_{x_{i} \in J(\theta)}
    \log
    f(x_i ; \theta)
   -
   \frac{1}{n}R_n(J(\theta)) \log{G}
   \; . 
   \nonumber
   \label{loglikelihood:proof:eq}
\end{eqnarray}
\qed

We want to bound the terms on the right hand side of 
(\ref{loglikelihood:lem:eq}) from above.
The first term is  easy.  In fact 
by  (\ref{permissible_quantity:eq}) and the strong law of large
numbers we have
\begin{equation}
\label{bounding-barf:eq}
\lim_{n\rightarrow\infty}  \frac{1}{n}\sum_{i=1}^{n}
  \log{
  \left\{
     \barf(x_{i};\btheta^{(s)},\rho(\btheta^{(s)})) + G
  \right\}
  }
 = 
 E_0[\log f(x;\theta_0)] - \lambda, \quad a.e.
\end{equation}

Next we consider the third term.  We prove the following lemma.

\begin{figure}[ht]
 \begin{center}
  \input{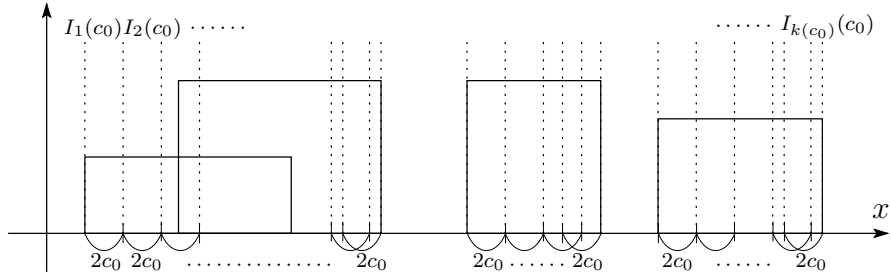}
 \end{center}
 \caption{Division of $J_{0}$ by short intervals of length $2c_{0}$.}
 \label{d_int:tpic:fig}
\end{figure}

\begin{lem}
\label{boundedRJ:lem}
 \begin{eqnarray}
 \limsup_{n \rightarrow \infty}
  \sup_{\theta \in \Theta_{n,K,s}'}
   \frac{1}{n} R_n(J(\theta))
   \leq 3M \cdot u \cdot 2c_{0}, 
  \quad a.e.
  \nonumber 
\end{eqnarray}
\end{lem}
\Proof Let $\epsilon > 0$ be arbitrarily fixed and let $J_{0}$ be the
support of the true density.  $J_{0}$ consists of at most $M$ intervals.  We
divide $J_{0}$ from $L_{\min}$ to $L_{\max}$ by short intervals of 
length $2c_{0}$.  In each right end of the intervals of $J_{0}$,
overlap of two short intervals of length $2c_{0}$ is allowed and the
right end of a short interval coincides with the right end of an
interval of $J_0$. See Figure \ref{d_int:tpic:fig}.
Let $k({c_{0}})$ be the number of
short intervals and let
$I_{1}(c_{0}),\ldots,I_{k(c_{0})}(c_{0})$ be the divided short
intervals.  Because $J_{0}$ consists of at most $M$ intervals, we
have
\begin{eqnarray}
 k(c_{0}) \leq 
   \frac{L}{2c_{0}} + M
  \; .
  \nonumber 
\end{eqnarray} 
Note that any interval in $J_{0}$ of length $2c_{0}$ is 
covered by at most $3$ small intervals 
from $\{I_{1}(c_{0}),\ldots,I_{k(c_{0})}(c_{0})\}$.
Now consider $J(\theta)$, the support of
$\sum_{m=1}^{K} \alpha_{m}f_{m}(x;\eta_{m})$.  The support of each
$f_{m}(x;\eta_{m})$, $1\le m \le K$, is an interval of length less
than or equal to  $2c_0$.  Therefore 
$J(\theta)$ is covered by at most $3M$ short intervals. 
Then 
the following relation holds.
\begin{eqnarray}
 \lefteqn{
  \sup_{\theta \in \Theta_{n,K,s}'}
  \frac{1}{n} R_n(J(\theta))
  - 3M \cdot u \cdot 2c_{0}
  > \epsilon
  } & &
  \\
 & \Rightarrow &
  {1} \leq \exists k \leq {k(c_{0})} \; , \;  
   \frac{1}{n} R_n(I_{k}(c_{0}))
   - u \cdot 2c_{0}
   > \frac{\epsilon}{3M}
  \; . 
  \label{supR_InequalityUsingSmallIntervals:eq}
  \nonumber
\end{eqnarray}
{}From (\ref{supR_InequalityUsingSmallIntervals:eq}), we have 
\begin{eqnarray}
 \lefteqn{
 \prob
  \left(
   \sup_{\theta \in \Theta_{n,K,s}'}
  \frac{1}{n} R_n(J(\theta))
  - 3M \cdot u \cdot 2c_{0}
  > \epsilon
  \right)
  } & & \nonumber \\
  & \leq &
   \sum_{k=1}^{k(c_{0})}
   \prob
   \left(
    \frac{1}{n} R_n(I_{k}(c_{0}))
    - u \cdot 2c_{0}
    > \frac{\epsilon}{3M}
   \right)
   \; . 
   \nonumber 
\end{eqnarray}
For any set $V\subset \real$, let $P_0(V)$ denote the probability of
$V$ under the true density
\begin{eqnarray}
 P_0(V) \eqv \int_{V} f(x;\theta_{0}) dx
  \; . 
  \nonumber
\end{eqnarray}
Then 
\begin{eqnarray}
 P_0(I_{k}(c_{0})) \leq u \cdot 2c_{0} , 
  \quad 
  k=1,\ldots,k(\theta)
  \; .
\end{eqnarray} 
Since $R_n(V) \sim \bin(n , P_0(V))$ 
and from (\ref{okamoto1:ineq}), 
we obtain \begin{eqnarray}
 \lefteqn{
  \prob
  \left(
   \frac{1}{n} R_n(I_{k}(c_{0}))
   - u \cdot 2c_{0}
   > \frac{\epsilon}{3M}
  \right)
  } & & \nonumber \\
 & \leq &
  \prob
  \left(
   \frac{1}{n} R_n(I_{k}(c_{0}))
   - P_0(I_{k}(c_{0}))
   > \frac{\epsilon}{3M}
  \right)
  \nonumber \\
 & \leq & 
  \exp{
  \left(
   -\frac{2n\epsilon^2}{9M^{2}}
  \right)} .
  \nonumber 
\end{eqnarray}
Therefore
\begin{eqnarray}
 {
  \prob
  \left(
   \sup_{\theta \in \Theta_{n,K,s}'}
   \frac{1}{n} R_n(J(\theta))
   - 3M \cdot u \cdot 2c_{0}
   > \epsilon
  \right)
  } 
  \leq 
  \left(
   \frac{L}{2c_{0}} + M
  \right)
  \exp{
  \left(
   -\frac{2n\epsilon^2}{9M^{2}}
  \right)}
  \; .
  \nonumber 
\end{eqnarray}
When we sum this over $n$, 
the resulting series on the right converges.
Hence by Borel-Cantelli, we have 
\begin{eqnarray}
 \prob
  \left(
    \sup_{\theta \in \Theta_{n,K,s}'}
    \frac{1}{n} R_n(J(\theta))
    - 3M \cdot u \cdot 2c_{0}
   > \epsilon
   \quad i.o.
  \right) = 0.
  \nonumber 
\end{eqnarray}
Because $\epsilon > 0$ was arbitrary, 
we obtain 
\begin{eqnarray}
 \limsup_{n \rightarrow \infty}
 \sup_{\theta \in \Theta_{n,K,s}'}
  \frac{1}{n} R_n(J(\theta))
  \leq
  3M \cdot u \cdot 2c_{0}, 
  \quad a.e. 
  \nonumber
\end{eqnarray}
\qed

By this lemma 
and
(\ref{1:c_0:condition:proof:eq}) we have
\begin{equation}
\label{boundingRJ:eq}
\limsup_{n \rightarrow \infty}
  \sup_{\theta \in \Theta_{n,K,s}'}
   \frac{1}{n} R_n(J(\theta)) \cdot ( - \log G )
   \leq 3M \cdot u \cdot 2c_{0}\cdot (-\log G)  < \frac{\lambda}{4}.
\end{equation}
This bounds the third term on the right hand side of 
(\ref{loglikelihood:lem:eq}) from above.

Finally we bound the second term on the right hand side of
(\ref{loglikelihood:lem:eq}) from above.  This is the most difficult
part of our proof.
For $x \in J(\theta)$ 
write $f(x;\theta)=\sum_{m=1}^{M} \alpha_m f_{m}(x;\eta_{m})$ as 
\begin{equation}
\label{eq:f-with-HJt}
f(x; \theta) = 
  \frac{1}{n}
  \sum_{t=1}^{T(\theta)}H(J_{t}(\theta))\indicator_{J_{t}(\theta)}(x) ,
\end{equation}
where $J_{t} \eqv J_{t}(\theta)$ are disjoint half-open intervals, 
$\indicator_{J_{t}(\theta)}(x)$ is the indicator function, 
$$
H(J_{t}(\theta))= f(x;\theta), \quad x\in J_t(\theta),
$$
is the height of $f(x;\theta)$ on $J_{t}(\theta)$ and 
$T \eqv T(\theta)$ is the number of the intervals $J_{t}(\theta)$.
Note that $T(\theta)\le 2M$, because
$f(x;\theta)$ changes its height only at $a_m - b_m$ 
or $a_m + b_m$, $m=1,\ldots,M$. 
For convenience we determine the order of $t$ such that    
\begin{eqnarray}
 H(J_{1}(\theta)) \leq H(J_{2}(\theta)) \leq
  \cdots \leq H(J_{T(\theta)}(\theta))
  \nonumber 
  \; . 
\end{eqnarray}
We now classify the intervals $J_{t}(\theta),\ t=1,\ldots, T(\theta),$
by the height $H(J_t(\theta))$.
Define $c_{n}'$ by
\begin{eqnarray}
 c_{n}' = c_{0} \cdot \exp{(-n^{1/4})}
  \nonumber
\end{eqnarray}
and define $\tau_n(\theta)$
\begin{eqnarray}
 \tau_n(\theta) \eqv
 \max \{
 t \in \{1,\ldots,T\}
 \mid
 H(J_{t}(\theta)) \leq \frac{M}{2c_{n}'}
 \} .
\end{eqnarray}
Then the second term  
on the right hand side of
(\ref{loglikelihood:lem:eq}) is 
written as
\begin{eqnarray}
 \frac{1}{n} \sum_{x_{i} \in J(\theta)}
  \log
  f(x_i ;\theta)
& = &
   \sum_{t=1}^{T(\theta)}
   \frac{1}{n}
   \sum_{x_{i} \in J_{t}(\theta)}
   \log{ H(J_{t}(\theta)) }
\label{scond-term-logH:eq} \\
& =&
  \frac{1}{n} 
  \sum_{t=1}^{T(\theta)} 
  R_n({J_{t}(\theta)}) \cdot \log{H(J_{t}(\theta))} 
  \nonumber \\
  & = &
  \frac{1}{n} 
  \sum_{t=1}^{\tau_{n}(\theta)} 
  R_n({J_{t}(\theta)}) \cdot \log{H(J_{t}(\theta))} 
 \nonumber \\
&& \qquad  
+ \frac{1}{n} \sum_{t=\tau_{n}(\theta) + 1}^{T(\theta)} 
  R_n({J_{t}(\theta)}) \cdot \log{H(J_{t}(\theta))} .
  \nonumber 
\end{eqnarray}
{}From (\ref{1:c_0:condition:proof:eq}),
(\ref{2:c_0:condition:proof:eq}), and noting that $\log x/x$ is
decreasing in $x\ge e$, 
we have 
\begin{eqnarray}
 3 \sum_{t=1}^{\tau_n(\theta)}
     \frac{u}{H(J_{t}(\theta))}
     \log{H(J_{t}(\theta))}
   & \leq &
    3 \cdot 2{M} \cdot u \cdot 2c_{0} \cdot \log{\frac{1}{2c_{0}}} 
   < 
  \frac{\lambda}{4} , 
  \nonumber \\
 \sum_{t=\tau_n(\theta)+1}^{T(\theta)}
  3 \cdot \frac{2}{n}\log{H(J_{t}(\theta))}
  &\leq&
  3 \cdot 2{M} \cdot \frac{2}{n} \cdot (n^d - \log\frac{M}{2c_0})
  \rightarrow 0 .
  \label{1:limit:ineq:eq}
\end{eqnarray}
Suppose that the following inequality holds.
\begin{eqnarray}
 \hspace{1cm}
 \lefteqn{
  \limsup_{n \rightarrow \infty}
  \sup_{\theta \in \Theta_{n,K,s}'}
  \left[
   \sum_{t=1}^{T(\theta)}
   \frac{1}{n}
   R_{n}(J_{t}(\theta))
   \log{ H(J_{t}(\theta)) }
   \right.
  } 
  & & 
  \\
 & & 
  - \left.
   3\left\{
     \sum_{t=1}^{\tau_n(\theta)}
     \frac{u}{H(J_{t}(\theta))}
     \log{H(J_{t}(\theta))}
     +
     \sum_{t=\tau_n(\theta)+1}^{T(\theta)}
     \frac{2}{n}\log{H(J_{t}(\theta))}
    \right\}
  \right]
  \leq 0 , 
  \quad a.e.
  \nonumber 
  \label{2:goal:eq}
\end{eqnarray}
Then from (\ref{scond-term-logH:eq})
and (\ref{1:limit:ineq:eq}), 
the second term on the right hand side of 
(\ref{loglikelihood:lem:eq}) is bounded from above as
\begin{equation}
\label{eq:bound-2nd-term:eq}
\limsup_{n\rightarrow\infty}\frac{1}{n} \sup_{\theta\in
  \Theta_{n,K,s}'}
\sum_{x_i \in J(\theta)} \log
f(x_i; \theta) \le \frac{4}{\lambda}.
\end{equation}
Combining (\ref{bounding-barf:eq}), (\ref{boundingRJ:eq}) and
(\ref{eq:bound-2nd-term:eq}) we obtain
\begin{eqnarray}
  \limsup_{n \rightarrow \infty}
  \sup_{\theta \in \Theta_{n,K,s}'}
  \frac{1}{n}  \sum_{i=1}^{n} 
  \log{f(x_i; \theta)}
 & \leq &
  \left( 
   E_{0}[\log{f(x;\theta_0)}] - \lambda 
  \right)
  + 
  \frac{\lambda}{4}
  + 
  \frac{\lambda}{4}
  \nonumber \\
 & \leq &
E_{0}[\log{f(x;\theta_0)}]
   - \frac{\lambda}{2} ,
   \quad a.e.
   \nonumber 
\end{eqnarray}
and ~(\ref{goal:eq}) is satisfied.  Therefore it suffices to
prove (\ref{2:goal:eq}), which is a new goal of our proof.

We now consider further finite covering of $\Theta_{n,K,s}'$.
Define 
\begin{eqnarray}
 \Theta_{n,K,s,T,\tau}'
  \eqv
  \{
  \theta \in \Theta_{n,K,s}'
  \mid 
  T(\theta) = T \; , \; \tau_{n}(\Theta) = \tau 
  \}
 \; .
 \nonumber 
\end{eqnarray}
Then 
\begin{eqnarray}
\label{bounding-2nd-with-fixed-T-tau:eq}
 \hspace{1cm} 
 \lefteqn{
  \sup_{\theta \in \Theta_{n,K,s}'}
  \left[
   \sum_{t=1}^{T(\theta)}
   \frac{1}{n}
   R_{n}(J_{t}(\theta))
   \log{ H(J_{t}(\theta)) }
   \right.
  } & & 
  \\
 & & \quad
 -   \left.
   3\left\{
     \sum_{t=1}^{\tau_n(\theta)}
     \frac{u}{H(J_{t}(\theta))}
     \log{H(J_{t}(\theta))}
     +
     \sum_{t=\tau_n(\theta)+1}^{T(\theta)}
     \frac{2}{n}\log{H(J_{t}(\theta))}
    \right\}
  \right]
\nonumber
\\
 &\leq &
  \max_{T=1,\ldots,2M} \max_{\tau=1,\ldots,T}
  \Biggl[
     \nonumber \\
 & & 
  \left.   \sup_{\theta \in \Theta_{n,K,s,T,\tau}'}
   \left\{
    \sum_{t=1}^{\tau}
    \frac{1}{n}
    R_{n}(J_{t}(\theta))
    \log{ H(J_{t}(\theta)) }
  - 
    3 \sum_{t=1}^{\tau}
    \frac{u}{H(J_{t}(\theta))}
    \log{H(J_{t}(\theta))}
   \right\} \right.
  \nonumber \\
 & & + 
  \sup_{\theta \in \Theta_{n,K,s,T,\tau}'}
  \left.
   \left\{
    \sum_{t=\tau + 1}^{T}
    \frac{1}{n}
    R_{n}(J_{t}(\theta))
    \log{ H(J_{t}(\theta)) }
    -
    3 \sum_{t=\tau + 1}^{T}
    \frac{2}{n}\log{H(J_{t}(\theta))}
   \right\}
  \right]
  \; . 
  \nonumber 
\end{eqnarray}
Suppose that  the following inequalities hold for all $T$ 
and $\tau$.
\begin{eqnarray}
\label{99-1:goal:eq}
 \lefteqn{
  \limsup_{n \rightarrow \infty}
  \sup_{\theta \in \Theta_{n,K,s,T,\tau}'}
  \left[
   \sum_{t=1}^{\tau}
   \frac{1}{n}
   R_{n}(J_{t}(\theta))
   \log{ H(J_{t}(\theta)) }
  \right. 
  }& &
  \\
 & & \hspace{4cm}
  \left.
   -
   3 \sum_{t=1}^{\tau}
   \frac{u}{H(J_{t}(\theta))}
   \log{H(J_{t}(\theta))}
  \right]
  \leq 0 , 
  \quad a.e.
  \nonumber 
\end{eqnarray}
\begin{eqnarray}
 \limsup_{n \rightarrow \infty}
  \sup_{\theta \in \Theta_{n,K,s,T,\tau}'}
  \left[
   \sum_{t=\tau + 1}^{T}
   \frac{1}{n}
   R_{n}(J_{t}(\theta))
   \log{ H(J_{t}(\theta)) }
   -
   3 \sum_{t=\tau + 1}^{T}
   \frac{2}{n}\log{H(J_{t}(\theta))}
  \right]
  \leq 0 , 
  \quad a.e.
  \nonumber \\
  \label{99-2:goal:eq}
\end{eqnarray}
Then (\ref{2:goal:eq}) is derived from 
(\ref{bounding-2nd-with-fixed-T-tau:eq}), 
(\ref{99-1:goal:eq}), 
(\ref{99-2:goal:eq}).
Therefore it suffices to prove (\ref{99-1:goal:eq}) and 
(\ref{99-2:goal:eq}), which are the final goals of our 
proof.  
We state (\ref{99-1:goal:eq}) and (\ref{99-2:goal:eq}) as 
two lemmas
and give their proofs.
\begin{lem}
 \begin{eqnarray}
  \limsup_{n \rightarrow \infty}
   \sup_{\theta \in \Theta_{n,K,s,T,\tau}'}
   \left[
    \sum_{t=\tau + 1}^{T}
    \frac{1}{n}
    R_{n}(J_{t}(\theta))
    \log{ H(J_{t}(\theta)) }
    -
    3 \sum_{t=\tau + 1}^{T}
    \frac{2}{n}\log{H(J_{t}(\theta))}
  \right]
   \leq 0
   \quad a.e.
   \nonumber 
  \label{99-2:goal:lem:eq}
 \end{eqnarray}
 \label{99-2:goal:lem}
\end{lem}

\Proof
Let $\delta > 0$ be  any fixed positive real constant and  
let $a_{t}'(\theta)$ denote the middle point of $J_{t}(\theta)$. 
Here, we consider the probability of the event that 
\begin{eqnarray}
 \hspace{1cm}
 \sup_{\theta \in \Theta_{n,K,s,T,\tau}'}
  \left[
   \sum_{t=\tau + 1}^{T}
   \frac{1}{n}
   R_{n}(J_{t}(\theta))
   \log{ H(J_{t}(\theta)) }
   -
   3 \sum_{t=\tau + 1}^{T}
   \frac{2}{n}\log{H(J_{t}(\theta))}
  \right]
  > 2M\delta .
  \label{99-2:RareEvent:eq}
\end{eqnarray}
Noting that for $t>\tau$, the length of $J_t(\theta)$ is less than or
equal to $2c_n'$, the following relation holds for this event.
\begin{eqnarray}
 \lefteqn{
  \textrm{The event (\ref{99-2:RareEvent:eq}) occurs. }
  } & & 
  \nonumber \\
 & \Rightarrow &
 \sup_{\theta \in \Theta_{n,K,s,T,\tau}'}
  \left[
   \sum_{t=\tau+1}^T
   \max
   \left\{
    0,
    \left(
     \frac{1}{n}
      R_n([a_{t}'(\theta) - c_{n}', a_{t}'(\theta) + c_{n}'])
    \right.
   \right.
  \right.
  \nonumber \\
 & & \hspace{5cm}
  \left.
   \left.
    \left.
     - 3 \cdot \frac{2}{n}
    \right)
   \right\}
   \log \frac{M}{2c_n}
  \right] 
 > 2M\delta
 \nonumber \\
 & \Rightarrow & 
  \exists \theta \in \Theta_{n,K,s,T,\tau}' ,
  \exists t > \tau 
  \nonumber \\ 
 & &
  \max
  \left\{
   0,
   \left(
    \frac{1}{n}
    R_n([a_{t}'(\theta) - c_{n}', a_{t}'(\theta) + c_{n}'])
    - 3 \cdot \frac{2}{n}
   \right)
  \right\}
  \log \frac{M}{2c_n}
  > \delta 
  \nonumber \\
 & \Rightarrow &
  \exists \theta \in \Theta_{n,K,s,T,\tau}' ,
  \exists t > \tau 
  \nonumber \\ 
 & &
  R_n([a_{t}'(\theta) - c_{n}', a_{t}'(\theta) + c_{n}'])
  \geq {6}
  \nonumber \\ 
 & \Rightarrow &
  \sup_{L_{\min} \leq a' \leq L_{\max}} R_n([a' - c_{n}', a' + c_{n}']) \geq 6
   \; .
  \label{99-2:rarer-event:eq}
\end{eqnarray}
Below, we consider the probability of the event 
that (\ref{99-2:rarer-event:eq}) occurs. 
We divide $J_{0}$ from $L_{\min}$ to $L_{\max}$ 
by short intervals of length $2c_{n}'$ as 
in the proof of Lemma \ref{boundedRJ:lem}.  
Let $k({c_{n}'})$ be the number of
short intervals and let
$I_{1}(c_{n}'),\ldots,I_{k(c_{n}')}(c_{n}')$ be the divided short
intervals. 
Because $J_{0}$ consists of at most $M$ intervals, we
have 
\begin{eqnarray}
 k(c_{n}') \leq 
   \frac{L}{2c_{n}'} + M
  \; .
  \label{99-2:bound-k_c_n:eq}
\end{eqnarray} 
Since any interval in $J_{0}$ of length $2c_{n}'$ is 
covered by at most $3$ small intervals 
from $\{I_{1}(c_{n}'),\ldots,I_{k(c_{n}')}(c_{n}')\}$
, the following relation holds. 
\begin{eqnarray}
 \hspace{1cm}
 \sup_{L_{\min} \leq a' \leq L_{\max}} 
  R_n([a' - c_{n}', a' + c_{n}']) \geq 6
  \Rightarrow
  {1} \leq \exists k \leq {k(c_{n}')} \; ,\; 
   R_n(I_{k}(c_{n}')) \geq 2
  \; . 
  \label{99-2:rarest-event:eq}
\end{eqnarray}
Note that $R_n(I_{k}(c_{n}'))
\sim \bin(n , P_0(I_{k}(c_{n}')))$ and 
$P_0(I_{k}(c_{n}')) \le 2c_n' u$. 
Therefore from 
 (\ref{99-2:rarer-event:eq}), (\ref{99-2:bound-k_c_n:eq}) 
 and~(\ref{99-2:rarest-event:eq}) 
we have 
\begin{eqnarray}
 \lefteqn{\hspace{-1cm}
  \prob
  \left(
   \sup_{\theta \in \Theta_{n,K,s,T,\tau}'}
   \left\{
    \sum_{t=\tau + 1}^{T}
    \frac{1}{n}
    R_{n}(J_{t}(\theta))
    \log{ H(J_{t}(\theta)) }
    -
    3 \sum_{t=\tau + 1}^{T}
    \frac{2}{n}\log{H(J_{t}(\theta))}
   \right\}
   > 2M\delta
  \right)
  } & & \nonumber \\
 &  & \hspace{4cm} \leq
  \left(
   \frac{L}{2c_{n}'} + M
  \right)
  \sum_{k=2}^{n}
  \begin{pmatrix} n \\ k \end{pmatrix}
  (2c_{n}' u)^{k}(1 - 2c_{n}' u)^{n-k}
   \nonumber \\
 &  & \hspace{4cm} \leq
  \left(
   \frac{L}{2c_{n}'} + M
  \right)
  \sum_{k=2}^{n}\frac{n^{k}}{k!}(2c_{n}'u)^{k}
  \nonumber \\
 &  & \hspace{4cm} \leq
    \left(
   \frac{L}{2c_{n}'} + M
  \right)
  (2nc_{n}'u)^{2}
  \exp{(2nc_{n}'u)}
  \; .
  \label{99-2:bounded-rare-event-prob:eq}
  \nonumber 
\end{eqnarray}
When we sum this over $n$, resulting series on the right converges.
Hence by Borel-Cantelli and the fact that $\delta > 0$ was arbitrary, 
we obtain 
\begin{eqnarray}
  \limsup_{n \rightarrow \infty}
   \sup_{\theta \in \Theta_{n,K,s,T,\tau}'}
   \left[
    \sum_{t=\tau + 1}^{T}
    \frac{1}{n}
    R_{n}(J_{t}(\theta))
    \log{ H(J_{t}(\theta)) }
    -
    3 \sum_{t=\tau + 1}^{T}
    \frac{2}{n}\log{H(J_{t}(\theta))}
  \right]
   \leq 0
   \quad a.e.
   \nonumber 
 \end{eqnarray}
\qed 

Finally we prove (\ref{99-1:goal:eq}).
\begin{lem}
 \begin{equation}
  \limsup_{n \rightarrow \infty}
  \sup_{\theta \in \Theta_{n,K,s,T,\tau}'}
  \left[
   \sum_{t=1}^{\tau}
   \frac{1}{n}
   R_{n}(J_{t}(\theta))
   \log{ H(J_{t}(\theta)) }
   -
   3 \sum_{t=1}^{\tau}
   \frac{u}{H(J_{t}(\theta))}
   \log{H(J_{t}(\theta))}
  \right]
  \leq 0
  \quad a.e.
  \nonumber
 \label{99-1:goal:lem:eq}
 \end{equation}
 \label{99-1:goal:lem}
\end{lem}

\Proof 
Let $\delta > 0$ be any fixed positive real constant 
and let $h_{n}$ be 
\begin{eqnarray}
 h_{n}
  \eqv
  \frac{\delta}{12}
  \left\{
   u \log{
   \left(
    \frac{M}{c_{n}'}
   \right)
   }
  \right\}^{-1}
 \; .
 \label{99-1:h_n:def:eq}
\end{eqnarray}
We divide $[c_{n}'/M, c_{0}]$ from $c_{0}$ to $c_{n}'/M$ 
by short intervals of length $h_{n}$.  
In the left end $c_{n}'/M$ of the interval $[c_{n}'/M, c_{0}]$,
overlap of two short intervals of length $h_{n}$ is allowed and the
left  end of a short interval is equal to $c_{n}'/M$.
Let $l_{n}$ be the number of
short intervals of length $h_{n}$ and define $b_{l}^{(n)}$ by
\begin{eqnarray}
 b_{l}^{(n)} \eqv 
 \begin{cases}
  c_{0} - (l - 1)h_{n}, & 1 \leq l \leq l_{n},\\
  {c_{n}'}/{M}, & l = l_{n} + 1.
 \end{cases}
 \nonumber
\end{eqnarray}
Then we have 
\begin{eqnarray}
 l_{n} \leq 
   \frac{c_{0}}{h_{n}} + 1
  \; . 
  \label{99-1:l_n:ineq:eq}
\end{eqnarray}
Next, we consider the probability of the event that 
\begin{eqnarray}
 \hspace{1cm}
   \sup_{\theta \in \Theta_{n,K,s,T,\tau}'}
  \left[
   \sum_{t=1}^{\tau}
   \frac{1}{n}
   R_{n}(J_{t}(\theta))
   \log{ H(J_{t}(\theta)) }
   -
   3 \sum_{t=1}^{\tau}
   \frac{u}{H(J_{t}(\theta))}
   \log{H(J_{t}(\theta))}
  \right]
  > 2M\delta .
  \label{99-1:rare-event:eq}
\end{eqnarray}
For this event the following relation holds.
\begin{eqnarray}
 \lefteqn{
  \textrm{The event (\ref{99-1:rare-event:eq}) occurs. }
  } & & 
  \nonumber \\
 & \Rightarrow &
  \exists \theta \in \Theta_{n,K,s,T,\tau}' 
  \; , \; 
  1 \leq 
  \exists l(1),\cdots, \exists l(\tau)
  \leq l_{n}
  \quad {\rm s.t.} \; 
  \nonumber \\
 & &
  2b_{l(1)+1}^{(n)} \leq \frac{1}{H(J_{1}(\theta))} \leq 2b_{l(1)}^{(n)},
  \cdots,
  2b_{l(\tau)+1}^{(n)} \leq \frac{1}{H(J_{\tau}(\theta))} \leq 2b_{l(\tau)}^{(n)},
  \nonumber \\
 & & 
   \sum_{t=1}^{\tau}
   \max
   \left\{
    0,
    \left(
     \frac{1}{n}
      R_n([a_{t}'(\theta) - b_{l(t)}^{(n)}, a_{t}'(\theta) + b_{l(t)}^{(n)}])
     -
     3u \cdot 2b_{l(t)+1}^{(n)}
    \right)
   \right\}
   \log{\frac{1}{2b_{l(t)+1}^{(n)}}}
   > 2M\delta
   \nonumber \\
 & \Rightarrow &
  \exists \theta \in \Theta_{n,K,s,T,\tau}' 
  \; , \; 
  1\leq \exists t \leq \tau   
  \; , \; 
  1 \leq 
  \exists l(t)
  \leq l_{n}
  \quad {\rm s.t.} \; 
  \nonumber \\
 & &
  2b_{l(t)+1}^{(n)} \leq \frac{1}{H(J_{t}(\theta))} \leq 2b_{l(t)}^{(n)},
  \nonumber \\
 & &
  \max
  \left\{
   0,
   \left(
    \frac{1}{n}
    R_n([a_{t}'(\theta) - b_{l(t)}^{(n)}, a_{t}'(\theta) + b_{l(t)}^{(n)}])
    -
    3u \cdot 2b_{l(t)+1}^{(n)}
   \right)
  \right\}
  \log{\frac{1}{2b_{l(t)+1}^{(n)}}}
  > \delta
  \nonumber \\
   & \Rightarrow &
  1 \leq 
  \exists l
  \leq l_{n}
  \quad {\rm s.t.} \; 
  \nonumber \\
 & &
  \max
  \left\{
   0,
   \sup_{L_{\min} \leq a' \leq L_{\max}}
   \left(
    \frac{1}{n}
    R_n([a' - b_{l}^{(n)}, a' + b_{l}^{(n)}])
    -
    3u \cdot 2b_{l+1}^{(n)}
   \right)
  \right\}
  \log{\frac{1}{2b_{l+1}^{(n)}}}
  > \delta
  \nonumber 
 \\
 & \Rightarrow &
  1 \leq \exists l \leq l_{n}
  \quad {\rm s.t.} \; 
  \nonumber \\
 & &
  \sup_{L_{\min} \leq a' \leq L_{\max}}
  \left\{
  \left(
   \frac{1}{n} R_n([a' - b_{l}^{(n)}, a' + b_{l}^{(n)}])
   -
   3u \cdot 2b_{l}^{(n)}
   \right)
   \log{\frac{1}{2b_{l+1}^{(n)}}}
  \right.
  \nonumber \\
 & &\qquad 
  +
  \left.
   3u(2b_{l}^{(n)} - 2b_{l+1}^{(n)})\log{\frac{1}{2b_{l+1}^{(n)}}}
  \right\}
  > \delta
  \nonumber \\
 \label{99-1:rarer-event:eq}
\end{eqnarray}
Then from (\ref{99-1:h_n:def:eq}) the following relation holds.
\begin{eqnarray}
 \lefteqn{
  \textrm{The event (\ref{99-1:rarer-event:eq}) occurs.}
  } & &
 \nonumber \\
 & \Rightarrow &
  {1}\leq \exists l \leq {l_{n}} \; , \;
   \sup_{L_{\min} \leq a' \leq L_{\max}}
   \frac{1}{n}
   \left(
    R_n([a' - b_{l}^{(n)},a' + b_{l}^{(n)}])
    -
    3u \cdot 2b_{l}^{(n)}
   \right)
   \log{\frac{1}{2b_{l+1}^{(n)}}}
   > \frac{\delta}{2}
  \nonumber \\
  \label{99-1:rarest-event:eq}
\end{eqnarray}

Below, we consider the probability of the event 
that (\ref{99-1:rarest-event:eq}) occurs.
We divide $J_{0}$ from $L_{\min}$ to $L_{\max}$ 
by short intervals of length $2b_{l}^{(n)}$ as 
in the proof of Lemma \ref{boundedRJ:lem}.  
Let $k({b_{l}^{(n)}})$ be the number of
short intervals and let
$I_{1}(b_{l}^{(n)}),\ldots,I_{k(b_{l}^{(n)})}(b_{l}^{(n)})$ be the divided short
intervals. 
Then we have 
\begin{eqnarray}
 k(b_l^{(n)}) \leq 
   \frac{L}{2b_{l}^{(n)}} + M
  \; .
  \label{99-1:bound-k_c_n:eq}
\end{eqnarray} 
Since any interval in $J_{0}$ of length $2b_l^{(n)}$ is 
covered by at most $3$ small intervals 
from $\{I_{1}(b_{l}^{(n)}),\ldots,I_{k(b_{l}^{(n)})}(b_{l}^{(n)})\}$,
 the following relation holds. 
\begin{eqnarray}
 \lefteqn{
  \sup_{L_{\min} \leq a' \leq L_{\max}}
  \left(
   \frac{1}{n}
   R_n([a' - b_{l}^{(n)},a' + b_{l}^{(n)}])
   -
   3u \cdot 2b_{l}^{(n)}
  \right)
  > 
  \frac{\delta}{2}
   \left(
    \log{\frac{1}{2b_{l+1}^{(n)}}}
   \right)^{-1}
   } & & \nonumber \\
 & \Rightarrow &
  \max_{k = 1,\ldots,k(b_{l}^{(n)})}
  \left(
   \frac{1}{n}
   R_n(I_k(b_{l}^{(n)}))
   -
   u \cdot 2b_{l}^{(n)}
  \right)
  > 
  \frac{1}{3} \cdot \frac{\delta}{2}
  \left(
   \log{\frac{1}{2b_{l+1}^{(n)}}}
  \right)^{-1}
  \; . 
  \label{99-1:short-interval-event:eq}
\end{eqnarray}
Note that $R_n(I_k(b_{l}^{(n)})) \sim \bin(n, P_0(I_k(b_{l}^{(n)})))$ 
and $P_0(I_k(b_{l}^{(n)})) \leq u \cdot 2b_{l}^{(n)}$.
Therefore from (\ref{okamoto1:ineq}) and (\ref{99-1:bound-k_c_n:eq}) we have  
\begin{eqnarray}
 \lefteqn{
  \prob
  \left(
   \max_{k = 1,\ldots,k(b_{l}^{(n)})}
   \frac{1}{n}
   \left( 
    R_n(I_k(b_{l}^{(n)}))
    -
    u \cdot 2b_{l}^{(n)}
   \right)
   > 
   \frac{1}{3} \cdot \frac{\delta}{2}
   \left(
    \log{\frac{1}{2b_{l+1}^{(n)}}}
   \right)^{-1}
  \right)
  } & & \nonumber \\
  & & \hspace{4cm} 
   \leq 
   \left(
    \frac{L}{2b_{l}^{(n)}} + M
   \right)
   \exp{
   \left\{
    -2n \cdot \frac{\delta^{2}}{36}
    \left(
     \log{\frac{1}{2b_{l+1}^{(n)}}}
    \right)^{-2}
   \right\}
   }
   \nonumber \\
 & & \hspace{4cm} 
   \leq 
   \left(
    \frac{L}{2c_{n}'} + M
   \right)
   \exp{
   \left\{
    -2n \cdot \frac{\delta^{2}}{36}
    \left(
     \log{\frac{1}{2c_{n}'}}
    \right)^{-2}
   \right\}
   }
   \label{99-1:c_n_dash:prob:eq}
   \; . 
\end{eqnarray}
{}From (\ref{99-1:l_n:ineq:eq}), 
(\ref{99-1:rarer-event:eq}), 
(\ref{99-1:rarest-event:eq}),
(\ref{99-1:short-interval-event:eq}),
(\ref{99-1:c_n_dash:prob:eq}), 
we obtain 
\begin{eqnarray}
& & \hspace{-5mm} 
 \prob
  \left(
   \sup_{\theta \in \Theta_{n,K,s,T,\tau}'}
   \left[
    \sum_{t=1}^{\tau}
    \frac{1}{n}
    R_{n}(J_{t}(\theta))
    \log{ H(J_{t}(\theta)) }
    -
    3 \sum_{t=1}^{\tau}
    \frac{u}{H(J_{t}(\theta))}
    \log{H(J_{t}(\theta))}
   \right]
   > 2M\delta
  \right)
  \nonumber \\
 &  & \hspace{4cm} \leq 
  \left(
   \frac{c_{0}}{h_{n}} + 1
  \right)
  \left(
   \frac{L}{2c_{n}'} + M
  \right)
  \exp{
  \left\{
   -2n \cdot \frac{\delta^{2}}{36}
   \left(
    \log{\frac{1}{2c_{n}'}}
   \right)^{-2}
  \right\}
  }
  \; . 
  \nonumber 
\end{eqnarray}
When we sum this over $n$, 
the resulting series on the right converges.
Hence by Borel-Cantelli and the fact 
that $\delta > 0$ is arbitrary, we have 
\begin{equation}
 \limsup_{n \rightarrow \infty}
  \sup_{\theta \in \Theta_{n,K,s,T,\tau}'}
  \left[
   \sum_{t=1}^{\tau}
   \frac{1}{n}
   R_{n}(J_{t}(\theta))
   \log{ H(J_{t}(\theta)) }
   -
   3 \sum_{t=1}^{\tau}
   \frac{u}{H(J_{t}(\theta))}
   \log{H(J_{t}(\theta))}
  \right]
  \leq 0
  \quad a.e.
  \nonumber
\end{equation}
\qed

This completes the proof of theorem \ref{thm:main}.

\end{document}